\documentclass[12pt]{article}
\usepackage{amsmath,amssymb,amsthm,latexsym,amscd}

\title{Closures and generating sets related \\
to combinations of structures\footnote{{\em Mathematics Subject Classification.}
03C30, 03C15, 03C50, 54A05.
\newline\indent \ \ \ The research is partially supported by
Committee of Science in Education and Science Ministry of the
Republic of Kazakhstan, Grant No. 0830/GF4. } }
\author{Sergey V.
Sudoplatov\footnote{sudoplat@math.nsc.ru}}
\date{}
\begin{document}
\maketitle

\begin{abstract}
We investigate closure operators and describe their properties for
$E$-combinations and $P$-combinations of structures and their
theories. We prove, for $E$-combinations, that the existence of a
minimal generating set of theories is equivalent to the existence
of the least generating set, and characterize syntactically and
semantically the property of the existence of the least generating
set. For the class of linearly ordered language uniform theories
we solve the problem of the existence of least generating set with
respect to $E$-combinations and characterize that existence in
terms of orders.

{\bf Key words:} $E$-combination, $P$-combination, closure
operator, generating set, language uniform theory.
\end{abstract}

\section{Introduction and preliminaries}

We continue to study structural properties of $E$-combin\-a\-tions
and $P$-combin\-a\-tions of structures and their theories
\cite{cs}.

In Section 2, using the $E$-operators and $P$-operators we
introduce topologies (related to topologies in \cite{BaPl}) and
investigate their properties.

In Section 3, we prove, for $E$-combinations, that the existence
of a minimal generating set of theories is equivalent to the
existence of the least generating set, and characterize
syntactically and semantically the property of the existence of
the least generating set.

In Section 4, for the class of linearly ordered language uniform
theories, we solve the problem of the existence of least
generating set with respect to $E$-combinations and characterize
that existence in terms of orders.

In Section 5 we describe some properties of $e$-spectra for
$E$-combinations of linearly ordered language uniform theories.

\medskip
Throughout the paper we use the following terminology in
\cite{cs}.

Let $P=(P_i)_{i\in I}$, be a family of nonempty unary predicates,
$(\mathcal{A}_i)_{i\in I}$ be a family of structures such that
$P_i$ is the universe of $\mathcal{A}_i$, $i\in I$, and the
symbols $P_i$ are disjoint with languages for the structures
$\mathcal{A}_j$, $j\in I$. The structure
$\mathcal{A}_P\rightleftharpoons\bigcup\limits_{i\in
I}\mathcal{A}_i$\index{$\mathcal{A}_P$} expanded by the predicates
$P_i$ is the {\em $P$-union}\index{$P$-union} of the structures
$\mathcal{A}_i$, and the operator mapping $(\mathcal{A}_i)_{i\in
I}$ to $\mathcal{A}_P$ is the {\em
$P$-operator}\index{$P$-operator}. The structure $\mathcal{A}_P$
is called the {\em $P$-combination}\index{$P$-combination} of the
structures $\mathcal{A}_i$ and denoted by ${\rm
Comb}_P(\mathcal{A}_i)_{i\in I}$\index{${\rm
Comb}_P(\mathcal{A}_i)_{i\in I}$} if
$\mathcal{A}_i=(\mathcal{A}_P\upharpoonright
A_i)\upharpoonright\Sigma(\mathcal{A}_i)$, $i\in I$. Structures
$\mathcal{A}'$, which are elementary equivalent to ${\rm
Comb}_P(\mathcal{A}_i)_{i\in I}$, will be also considered as
$P$-combinations.

Clearly, all structures $\mathcal{A}'\equiv {\rm
Comb}_P(\mathcal{A}_i)_{i\in I}$ are represented as unions of
their restrictions $\mathcal{A}'_i=(\mathcal{A}'\upharpoonright
P_i)\upharpoonright\Sigma(\mathcal{A}_i)$ if and only if the set
$p_\infty(x)=\{\neg P_i(x)\mid i\in I\}$ is inconsistent. If
$\mathcal{A}'\ne{\rm Comb}_P(\mathcal{A}'_i)_{i\in I}$, we write
$\mathcal{A}'={\rm Comb}_P(\mathcal{A}'_i)_{i\in
I\cup\{\infty\}}$, where
$\mathcal{A}'_\infty=\mathcal{A}'\upharpoonright
\bigcap\limits_{i\in I}\overline{P_i}$, maybe applying
Morleyzation. Moreover, we write ${\rm
Comb}_P(\mathcal{A}_i)_{i\in I\cup\{\infty\}}$\index{${\rm
Comb}_P(\mathcal{A}_i)_{i\in I\cup\{\infty\}}$} for ${\rm
Comb}_P(\mathcal{A}_i)_{i\in I}$ with the empty structure
$\mathcal{A}_\infty$.

Note that if all predicates $P_i$ are disjoint, a structure
$\mathcal{A}_P$ is a $P$-combination and a disjoint union of
structures $\mathcal{A}_i$. In this case the $P$-combination
$\mathcal{A}_P$ is called {\em
disjoint}.\index{$P$-combination!disjoint} Clearly, for any
disjoint $P$-combination $\mathcal{A}_P$, ${\rm
Th}(\mathcal{A}_P)={\rm Th}(\mathcal{A}'_P)$, where
$\mathcal{A}'_P$ is obtained from $\mathcal{A}_P$ replacing
$\mathcal{A}_i$ by pairwise disjoint
$\mathcal{A}'_i\equiv\mathcal{A}_i$, $i\in I$. Thus, in this case,
similar to structures the $P$-operator works for the theories
$T_i={\rm Th}(\mathcal{A}_i)$ producing the theory $T_P={\rm
Th}(\mathcal{A}_P)$\index{$T_P$}, which is denoted by ${\rm
Comb}_P(T_i)_{i\in I}$.\index{${\rm Comb}_P(T_i)_{i\in I}$}

For an equivalence relation $E$ replacing disjoint predicates
$P_i$ by $E$-classes we get the structure
$\mathcal{A}_E$\index{$\mathcal{A}_E$} being the {\em
$E$-union}\index{$E$-union} of the structures $\mathcal{A}_i$. In
this case the operator mapping $(\mathcal{A}_i)_{i\in I}$ to
$\mathcal{A}_E$ is the {\em $E$-operator}\index{$E$-operator}. The
structure $\mathcal{A}_E$ is also called the {\em
$E$-combination}\index{$E$-combination} of the structures
$\mathcal{A}_i$ and denoted by ${\rm Comb}_E(\mathcal{A}_i)_{i\in
I}$\index{${\rm Comb}_E(\mathcal{A}_i)_{i\in I}$}; here
$\mathcal{A}_i=(\mathcal{A}_E\upharpoonright
A_i)\upharpoonright\Sigma(\mathcal{A}_i)$, $i\in I$. Similar
above, structures $\mathcal{A}'$, which are elementary equivalent
to $\mathcal{A}_E$, are denoted by ${\rm
Comb}_E(\mathcal{A}'_j)_{j\in J}$, where $\mathcal{A}'_j$ are
restrictions of $\mathcal{A}'$ to its $E$-classes.

Clearly, $\mathcal{A}'\equiv\mathcal{A}_P$ realizing $p_\infty(x)$
is not elementary embeddable into $\mathcal{A}_P$ and can not be
represented as a disjoint $P$-combination of
$\mathcal{A}'_i\equiv\mathcal{A}_i$, $i\in I$. At the same time,
there are $E$-combinations such that all
$\mathcal{A}'\equiv\mathcal{A}_E$ can be represented as
$E$-combinations of some $\mathcal{A}'_j\equiv\mathcal{A}_i$. We
call this representability of $\mathcal{A}'$ to be the {\em
$E$-representability}.

If there is $\mathcal{A}'\equiv\mathcal{A}_E$ which is not
$E$-representable, we have the $E'$-representability replacing $E$
by $E'$ such that $E'$ is obtained from $E$ adding equivalence
classes with models for all theories $T$, where $T$ is a theory of
a restriction $\mathcal{B}$ of a structure
$\mathcal{A}'\equiv\mathcal{A}_E$ to some $E$-class and
$\mathcal{B}$ is not elementary equivalent to the structures
$\mathcal{A}_i$. The resulting structure $\mathcal{A}_{E'}$ (with
the $E'$-representability) is a {\em
$e$-completion}\index{$e$-completion}, or a {\em
$e$-saturation}\index{$e$-saturation}, of $\mathcal{A}_{E}$. The
structure $\mathcal{A}_{E'}$ itself is called {\em
$e$-complete}\index{Structure!$e$-complete}, or {\em
$e$-saturated}\index{Structure!$e$-saturated}, or {\em
$e$-universal}\index{Structure!$e$-universal}, or {\em
$e$-largest}\index{Structure!$e$-largest}.

For a structure $\mathcal{A}_E$ the number of {\em
new}\index{Structure!new} structures with respect to the
structures $\mathcal{A}_i$, i.~e., of the structures $\mathcal{B}$
which are pairwise elementary non-equivalent and elementary
non-equivalent to the structures $\mathcal{A}_i$, is called the
{\em $e$-spectrum}\index{$e$-spectrum} of $\mathcal{A}_E$ and
denoted by $e$-${\rm Sp}(\mathcal{A}_E)$.\index{$e$-${\rm
Sp}(\mathcal{A}_E)$} The value ${\rm sup}\{e$-${\rm
Sp}(\mathcal{A}'))\mid\mathcal{A}'\equiv\mathcal{A}_E\}$ is called
the {\em $e$-spectrum}\index{$e$-spectrum} of the theory ${\rm
Th}(\mathcal{A}_E)$ and denoted by $e$-${\rm Sp}({\rm
Th}(\mathcal{A}_E))$.\index{$e$-${\rm Sp}({\rm
Th}(\mathcal{A}_E))$}

If $\mathcal{A}_E$ does not have $E$-classes $\mathcal{A}_i$,
which can be removed, with all $E$-classes
$\mathcal{A}_j\equiv\mathcal{A}_i$, preserving the theory ${\rm
Th}(\mathcal{A}_E)$, then $\mathcal{A}_E$ is called {\em
$e$-prime}\index{Structure!$e$-prime}, or {\em
$e$-minimal}\index{Structure!$e$-minimal}.

For a structure $\mathcal{A}'\equiv\mathcal{A}_E$ we denote by
${\rm TH}(\mathcal{A}')$ the set of all theories ${\rm
Th}(\mathcal{A}_i)$\index{${\rm Th}(\mathcal{A}_i)$} of
$E$-classes $\mathcal{A}_i$ in $\mathcal{A}'$.

By the definition, an $e$-minimal structure $\mathcal{A}'$
consists of $E$-classes with a minimal set ${\rm
TH}(\mathcal{A}')$. If ${\rm TH}(\mathcal{A}')$ is the least for
models of ${\rm Th}(\mathcal{A}')$ then $\mathcal{A}'$ is called
{\em $e$-least}.\index{Structure!$e$-least}

\section{Closure operators}

{\bf Definition.} Let $\overline{\mathcal{T}}$ be the class of all
complete elementary theories of relational languages. For a set
$\mathcal{T}\subset\overline{\mathcal{T}}$ we denote by ${\rm
Cl}_E(\mathcal{T})$ the set of all theories ${\rm
Th}(\mathcal{A})$, where $\mathcal{A}$ is a structure of some
$E$-class in $\mathcal{A}'\equiv\mathcal{A}_E$,
$\mathcal{A}_E={\rm Comb}_E(\mathcal{A}_i)_{i\in I}$, ${\rm
Th}(\mathcal{A}_i)\in\mathcal{T}$. As usual, if $\mathcal{T}={\rm
Cl}_E(\mathcal{T})$ then $\mathcal{T}$ is said to be {\em
$E$-closed}.\index{Set!$E$-closed}

\medskip
By the definition,
\begin{equation}\label{cs2}
{\rm Cl}_E(\mathcal{T})={\rm TH}(\mathcal{A}'_{E'}),
\end{equation}
where $\mathcal{A}'_{E'}$ is an $e$-largest model of ${\rm
Th}(\mathcal{A}_{E})$, $\mathcal{A}_{E}$ consists of $E$-classes
representing models of all theories in $\mathcal{T}$.

Note that the equality (\ref{cs2}) does not depend on the choice
of $e$-largest model of ${\rm Th}(\mathcal{A}_{E})$.

The following proposition is obvious.

\medskip
{\bf Proposition 2.1.} $(1)$ {\em If $\mathcal{T}_0$,
$\mathcal{T}_1$ are sets of theories,
$\mathcal{T}_0\subseteq\mathcal{T}_1\subset\overline{\mathcal{T}}$,
then $\mathcal{T}_0\subseteq{\rm Cl}_E(\mathcal{T}_0)\subseteq{\rm
Cl}_E(\mathcal{T}_1)$.

\medskip
$(2)$ For any set $\mathcal{T}\subset\overline{\mathcal{T}}$,
$\mathcal{T}\subset{\rm Cl}_E(\mathcal{T})$ if and only if the
structure composed by $E$-classes of models of theories in
$\mathcal{T}$ is not $e$-largest.

\medskip
$(3)$ Every finite set $\mathcal{T}\subset\overline{\mathcal{T}}$
is $E$-closed.

\medskip
$(4)$ {\rm (Negation of finite character)} For any $T\in{\rm
Cl}_E(\mathcal{T})\setminus\mathcal{T}$ there are no finite
$\mathcal{T}_0\subset\mathcal{T}$ such that $T\in{\rm
Cl}_E(\mathcal{T}_0).$

\medskip
$(5)$ Any intersection of $E$-closed sets is $E$-closed.}

\medskip
For a set $\mathcal{T}\subset\overline{\mathcal{T}}$ of theories
in a language $\Sigma$ and for a sentence $\varphi$ with
$\Sigma(\varphi)\subseteq\Sigma$ we denote by
$\mathcal{T}_\varphi$\index{$\mathcal{T}_\varphi$} the set
$\{T\in\mathcal{T}\mid\varphi\in T\}$. Denote by
$\mathcal{T}_F$\index{$\mathcal{T}_F$} the family of all sets
$\mathcal{T}_\varphi$.

Clearly, the partially ordered set
$\langle\mathcal{T}_F;\subseteq\rangle$ forms a Boolean algebra
with the least element $\varnothing=\mathcal{T}_{\neg(x\approx
x)}$, the greatest element $\mathcal{T}=\mathcal{T}_{(x\approx
x)}$, and operations $\wedge$, $\vee$, $\bar{\,\, }$ satisfying
the following equalities: $
\mathcal{T}_\varphi\wedge\mathcal{T}_\psi=\mathcal{T}_{(\varphi\wedge\psi)},$
$\mathcal{T}_\varphi\vee\mathcal{T}_\psi=\mathcal{T}_{(\varphi\vee\psi)},$
$\overline{\mathcal{T}_\varphi}=\mathcal{T}_{\neg\varphi}.$

\medskip
By the definition, $\mathcal{T}_\varphi\subseteq\mathcal{T}_\psi$
if and only if for any model $\mathcal{M}$ of a theory in
$\mathcal{T}$ satisfying $\varphi$ we have
$\mathcal{M}\models\psi$.

\medskip
{\bf Proposition 2.2.} {\em If
$\mathcal{T}\subset\overline{\mathcal{T}}$ is an infinite set and
$T\in\overline{\mathcal{T}}\setminus\mathcal{T}$ then $T\in{\rm
Cl}_E(\mathcal{T})$ {\rm (}i.e., $T$ is an {\sl accumulation
point} for $\mathcal{T}$ with respect to $E$-closure ${\rm
Cl}_E${\rm )} if and only if for any formula $\varphi\in T$ the
set $\mathcal{T}_\varphi$ is infinite.}

\medskip
{\bf\em Proof.} Assume that there is a formula $\varphi\in T$ such
that only finitely many theories in $\mathcal{T}$, say
$T_1,\ldots,T_n$, satisfy $\varphi$. Since $T\notin\mathcal{T}$
then there is $\psi\in T$ such that $\psi\notin T_1\cup\ldots\cup
T_n$. Then $(\varphi\wedge\psi)\in T$ does not belong to all
theories in $\mathcal{T}$. Since $(\varphi\wedge\psi)$ does not
satisfy $E$-classes in models of $T_E={\rm Comb}_E(T_i)_{i\in I}$,
where $\mathcal{T}=\{T_i\mid i\in I\}$, we have $T\notin{\rm
Cl}_E(\mathcal{T})$.

If for any formula $\varphi\in T$, $\mathcal{T}_\varphi$ is
infinite then $\{\varphi^E\mid\varphi\in T\}\cup T_E$ (where
$\varphi^E$ are $E$-relativizations of the formulas $\varphi$) is
locally satisfied and so satisfied. Since $T_E$ is a complete
theory then $\{\varphi^E\mid\varphi\in T\}\subset T_E$ and hence
$T\in{\rm Cl}_E(\mathcal{T})$.~$\Box$

\medskip
Proposition 2.2 shows that the closure ${\rm Cl}_E$ corresponds to
the closure with respect to the ultraproduct operator \cite{KeCh,
ErPa, Bank1, Bank2}.

\medskip
{\bf Theorem 2.3.} {\em For any sets
$\mathcal{T}_0,\mathcal{T}_1\subset\overline{\mathcal{T}}$, ${\rm
Cl}_E(\mathcal{T}_0\cup\mathcal{T}_1)={\rm
Cl}_E(\mathcal{T}_0)\cup{\rm Cl}_E(\mathcal{T}_1)$.}

\medskip
{\bf\em Proof.} We have ${\rm Cl}_E(\mathcal{T}_0)\cup{\rm
Cl}_E(\mathcal{T}_1)\subseteq {\rm
Cl}_E(\mathcal{T}_0\cup\mathcal{T}_1)$ by Proposition 2.1 (1).

Let $T\in {\rm Cl}_E(\mathcal{T}_0\cup\mathcal{T}_1)$ and we argue
to show that $T\in{\rm Cl}_E(\mathcal{T}_0)\cup{\rm
Cl}_E(\mathcal{T}_1)$. Without loss of generality we assume that
$T\notin\mathcal{T}_0\cup\mathcal{T}_1$ and by Proposition 2.1
(3), $\mathcal{T}_0\cup\mathcal{T}_1$ is infinite. Define a
function $f\mbox{: }T\to\mathcal{P}(\{0,1\})$ by the following
rule: $f(\varphi)$ is the set of indexes $k\in\{0,1\}$ such that
$\varphi$ belongs to infinitely many theories in $\mathcal{T}_k$.
Note that $f(\varphi)$ is always nonempty since by Proposition
2.2, $\varphi$ belong to infinitely many theories in
$\mathcal{T}_0\cup\mathcal{T}_1$ and so to infinitely many
theories in $\mathcal{T}_0$ or to infinitely many theories in
$\mathcal{T}_1$. Again by Proposition 2.2 we have to prove that
$0\in f(\varphi)$ for each formula $\varphi\in T$ or $1\in
f(\varphi)$ for each formula $\varphi\in T$. Assuming on contrary,
there are formulas $\varphi,\psi\in T$ such that
$f(\varphi)=\{0\}$ and $f(\psi)=\{1\}$. Since
$(\varphi\wedge\psi)\in T$ and $f(\varphi\wedge\psi)$ is nonempty
we have $0\in f(\varphi\wedge\psi)$ or $1\in
f(\varphi\wedge\psi)$. In the first case, since
$\mathcal{T}_{\varphi\wedge\psi}\subseteq\mathcal{T}_\psi$ we get
$0\in f(\psi)$. In the second case, since
$\mathcal{T}_{\varphi\wedge\psi}\subseteq\mathcal{T}_\varphi$ we
get $1\in f(\varphi)$. Both cases contradict the assumption. Thus,
$T\in{\rm Cl}_E(\mathcal{T}_0)\cup{\rm
Cl}_E(\mathcal{T}_1)$.~$\Box$

\medskip
{\bf Corollary 2.4.} (Exchange property) {\em If $T_1\in{\rm
Cl}_E(\mathcal{T}\cup\{T_2\})\setminus{\rm Cl}_E(\mathcal{T})$
then $T_2\in{\rm Cl}_E(\mathcal{T}\cup\{T_1\})$.}

\medskip
{\bf\em Proof.} Since $T_1\in{\rm
Cl}_E(\mathcal{T}\cup\{T_2\})={\rm Cl}_E(\mathcal{T})\cup\{T_2\}$
by Proposition 2.1 (3) and Theorem 2.3, and $T_1\notin{\rm
Cl}_E(\mathcal{T})$, then $T_1=T_2$ and $T_2\in{\rm
Cl}_E(\mathcal{T}\cup\{T_1\})$ in view of Proposition 2.1
(1).~$\Box$

\medskip
{\bf Definition} \cite{Engel}. A {\em topological
space}\index{Space topological} is a pair $(X,\mathcal{O})$
consisting of a set $X$ and a family $\mathcal{O}$ of {\em
open}\index{Subset open} subsets of $X$ satisfying the following
conditions:

(O1) $\varnothing\in\mathcal{O}$ and $X\in\mathcal{O}$;

(O2) If $U_1\in\mathcal{O}$ and $U_2\in\mathcal{O}$ then $U_1\cap
U_2\in\mathcal{O}$;

(O3) If $\mathcal{O}'\subseteq\mathcal{O}$ then
$\cup\mathcal{O}'\in\mathcal{O}$.

\medskip
{\bf Definition} \cite{Engel}. A topological space
$(X,\mathcal{O})$ is a {\em $T_0$-space}\index{$T_0$-space} if for
any pair of distinct elements $x_1,x_2\in X$ there is an open set
$U\in\mathcal{O}$ containing exactly one of these elements.

\medskip
{\bf Definition} \cite{Engel}. A topological space
$(X,\mathcal{O})$ is {\em Hausdorff}\index{Space
topological!Hausdorff} if for any pair of distinct points
$x_1,x_2\in X$ there are open sets $U_1,U_2\in\mathcal{O}$ such
that $x_1\in U_1$, $x_2\in U_2$, and $U_1\cap U_2=\varnothing$.

\medskip
Let $\mathcal{T}\subset\overline{\mathcal{T}}$ be a set,
$\mathcal{O}_E(\mathcal{T})=\{\mathcal{T}\setminus{\rm
Cl}_E(\mathcal{T}')\mid\mathcal{T}'\subseteq\mathcal{T}\}$.\index{$\mathcal{O}_E(\mathcal{T})$}
Proposition 2.1 and Theorem 2.3 imply that the axioms (O1)--(O3)
are satisfied. Moreover, since for any theory $T\in
\overline{\mathcal{T}}$, ${\rm Cl}_E(\{T\})=\{T\}$ and hence,
$\mathcal{T}\setminus{\rm Cl}_E(\{T\})=\mathcal{T}\{T\}$ is an
open set containing all theories in $\mathcal{T}$, which are not
equal to $T$, then $(\mathcal{T},\mathcal{O}_E(\mathcal{T}))$ is a
$T_0$-space. Moreover, it is Hausdorff. Indeed, taking two
distinct theories $T_1,T_2\in \mathcal{T}$ we have a formula
$\varphi$ such that $\varphi\in T_1$ and $\neg\varphi\in T_2$. By
Proposition 2.2 we have that $\mathcal{T}_\varphi$ and
$\mathcal{T}_{\neg\varphi}$ are closed containing $T_1$ and $T_2$
respectively; at the same time $\mathcal{T}_\varphi$ and
$\mathcal{T}_{\neg\varphi}$ form a partition of $\mathcal{T}$, so
$\mathcal{T}_\varphi$ and $\mathcal{T}_{\neg\varphi}$ are disjoint
open sets.  Thus we have

\medskip
{\bf Theorem 2.5.} {\em For any set
$\mathcal{T}\subset\overline{\mathcal{T}}$ the pair
$(\mathcal{T},\mathcal{O}_E(\mathcal{T}))$ is a Hausdorff
topological space.}

\medskip
Similarly to the operator ${\rm Cl}_E(\mathcal{T})$ we define the
operator ${\rm Cl}_P(\mathcal{T})$ for families $P$ of 
predicates $P_i$ as follows.

\medskip
{\bf Definition.} For a set
$\mathcal{T}\subset\overline{\mathcal{T}}$ we denote by ${\rm
Cl}_P(\mathcal{T})$\index{${\rm Cl}_P(\mathcal{T})$} the set of
all theories ${\rm Th}(\mathcal{A})$ such that ${\rm
Th}(\mathcal{A})\in\mathcal{T}$ or $\mathcal{A}$ is a structure of
type $p_\infty(x)$ in $\mathcal{A}'\equiv\mathcal{A}_P$, where
$\mathcal{A}_P={\rm Comb}_P(\mathcal{A}_i)_{i\in I}$ and ${\rm
Th}(\mathcal{A}_i)\in\mathcal{T}$ are pairwise distinct. As above,
if $\mathcal{T}={\rm Cl}_P(\mathcal{T})$ then $\mathcal{T}$ is
said to be {\em $P$-closed}.\index{Set!$P$-closed}

Using above only disjoint $P$-combinations $\mathcal{A}_P$ we get
the closure ${\rm Cl}^d_P(\mathcal{T})$\index{${\rm
Cl}^d_P(\mathcal{T})$} being a subset of ${\rm
Cl}_P(\mathcal{T})$.

\medskip
The following example illustrates the difference between ${\rm
Cl}_P(\mathcal{T})$ and ${\rm Cl}^d_P(\mathcal{T})$.

\medskip
{\bf Example 2.7.} Taking disjoint copies of predicates
$P_i=\{a\in M_0\mid a<c_i\}$ with their $<$-structures as in
\cite[Example 4.8]{cs}, ${\rm
Cl}^d_P(\mathcal{T})\setminus\mathcal{T}$ produces models of the
Ehrenfeucht example and unboundedly many connected components each
of which is a copy of a model of the Ehrenfeucht example. At the
same time ${\rm Cl}_P(\mathcal{T})$ produces two new structures:
densely ordered structures with and without the least element.

\medskip
The following proposition is obvious.

\medskip
{\bf Proposition 2.8.} $(1)$ {\em If $\mathcal{T}_0$,
$\mathcal{T}_1$ are sets of theories,
$\mathcal{T}_0\subseteq\mathcal{T}_1\subset\overline{\mathcal{T}}$,
then $\mathcal{T}_0\subseteq{\rm Cl}_P(\mathcal{T}_0)\subseteq{\rm
Cl}_P(\mathcal{T}_1)$.

\medskip
$(2)$ Every finite set $\mathcal{T}\subset\overline{\mathcal{T}}$
is $P$-closed.

\medskip
$(3)$ {\rm (Negation of finite character)} For any $T\in{\rm
Cl}_P(\mathcal{T})\setminus\mathcal{T}$ there are no finite
$\mathcal{T}_0\subset\mathcal{T}$ such that $T\in{\rm
Cl}_P(\mathcal{T}_0)$.

\medskip
$(4)$ Any intersection of $P$-closed sets is $P$-closed.}

\medskip
{\bf Remark 2.9.} Note that an analogue of Proposition 2.8 for
$P$-combinations fails. Indeed, taking disjoint predicates $P_i$,
$i\in\omega$, with $2i+1$ elements and with structures
$\mathcal{A}_i$ of the empty language, we get, for the set
$\mathcal{T}$ of theories $T_i={\rm Th}(\mathcal{A}_i)$, that
${\rm Cl}_P(\mathcal{T})$ consists of the theories whose models
have cardinalities witnessing all ordinals in $\omega+1$. Thus,
for instance, theories in $\mathcal{T}$ do not contain the formula
\begin{equation}\label{form_x}
\exists x,y(\neg(x\approx y)\wedge\forall z((z\approx
x)\vee(z\approx y)))
\end{equation}
whereas ${\rm Cl}_P(\mathcal{T})$ (which is equal to ${\rm
Cl}^d_P(\mathcal{T})$) contains a theory with the formula
(\ref{form_x}).

More generally, for ${\rm Cl}^d_P(\mathcal{T})$ with infinite
$\mathcal{T}$, we have the following.

Since there are no links between distinct $P_i$, the structures of
$p_\infty(x)$ are defined as disjoint unions of connected
components $C(a)$, for $a$ realizing $p_\infty(x)$, where each
$C(a)$ consists of a set of realizations of $p_\infty$-preserving
formulas $\psi(a,x)$ (i.e., of formulas $\varphi(a,x)$ with
$\psi(a,x)\vdash p_\infty(x)$). Similar to Proposition 2.2
theories $T_{\infty,C(a)}$ of $C(a)$-restrictions of
$\mathcal{A}_\infty$ coincide and are characterized by the
following property: $T_{\infty,C(a)}\in{\rm Cl}^d_P(\mathcal{T})$
if and only if $T_{\infty,C(a)}\in\mathcal{T}$ or for any formula
$\varphi\in T_{\infty,C(a)}$, there are infinitely many theories
$T$ in $\mathcal{T}$ such that $\varphi$ satisfies all structures
approximating $C(a)$-restrictions of models of $T$.

Thus similarly to 2.3--2.5 we get the following three assertions
for disjoint $P$-combinations.

\medskip
{\bf Theorem 2.10.} {\em For any sets
$\mathcal{T}_0,\mathcal{T}_1\subset\overline{\mathcal{T}}$, ${\rm
Cl}^d_P(\mathcal{T}_0\cup\mathcal{T}_1)={\rm
Cl}^d_P(\mathcal{T}_0)\cup{\rm Cl}^d_P(\mathcal{T}_1)$.}

\medskip
{\bf Corollary 2.11.} (Exchange property) {\em If $T_1\in{\rm
Cl}^d_P(\mathcal{T}\cup\{T_2\})\setminus{\rm Cl}^d_P(\mathcal{T})$
then $T_2\in{\rm Cl}^d_P(\mathcal{T}\cup\{T_1\})$.}

\medskip
Let $\mathcal{T}\subset\overline{\mathcal{T}}$ be a set,
$\mathcal{O}^d_P(\mathcal{T})=\{\mathcal{T}\setminus{\rm
Cl}^d_P(\mathcal{T}')\mid\mathcal{T}'\subseteq\mathcal{T}\}$.\index{$\mathcal{O}^d_P(\mathcal{T})$}

\medskip
{\bf Theorem 2.12.} {\em For any set
$\mathcal{T}\subset\overline{\mathcal{T}}$ the pair
$(\mathcal{T},\mathcal{O}^d_P(\mathcal{T}))$ is a topological
$T_0$-space.}

\medskip
{\bf Remark 2.13.} By Proposition 2.8 (2), for any finite
$\mathcal{T}$ the spaces
$(\mathcal{T},\mathcal{O}_P(\mathcal{T}))$ and
$(\mathcal{T},\mathcal{O}^d_P(\mathcal{T}))$ are Hausdorff,
moreover, here
$\mathcal{O}_P(\mathcal{T})=\mathcal{O}^d_P(\mathcal{T})$
consisting of all subsets of $\mathcal{T}$. However, in general,
the spaces $(\mathcal{T},\mathcal{O}_P(\mathcal{T}))$ and
$(\mathcal{T},\mathcal{O}^d_P(\mathcal{T}))$ are not Hausdorff.

Indeed, consider structures $\mathcal{A}_i$, $i\in I$, where
$I=(\omega+1)\setminus\{0\}$, of the empty language and such that
$|\mathcal{A}_i|=i$. Let $T_i={\rm Th}(\mathcal{A}_i)$, $i\in I$,
$\mathcal{T}=\{T_i\mid i\in I\}$. Coding the theories $T_i$ by
their indexes we have the following. For any finite set $F\subset
I$, ${\rm Cl}_P(F)={\rm Cl}^d_P(F)=F$, and for any infinite set
${\rm INF}\subseteq I$, ${\rm Cl}_P({\rm INF})={\rm Cl}^d_P({\rm
INF})=I$. So any open set $U$ is either cofinite or empty. Thus
any two nonempty open sets are not disjoint.

Notice that we get a similar effect replacing elements in
$\mathcal{A}_i$ by equivalence classes with pairwise isomorphic
finite structures, may be with additional classes having arbitrary
structures.

\medskip
{\bf Remark 2.14.} If the closure operator ${\rm
Cl}^{d,r}_P$\index{${\rm Cl}^{d,r}_P$} is obtained from ${\rm
Cl}^{d}_P$ permitting repetitions of structures for predicates
$P_i$, we can lose both the property of $T_0$-space and the
identical closure for finite sets of theories. Indeed, for the
example in Remark 2.13, ${\rm Cl}^{d,r}_P(\mathcal{T})$ is equal
to the ${\rm Cl}^{d,r}_P$-closure of any singleton $\{T\}\in {\rm
Cl}^{d,r}_P(\mathcal{T})$ since the type $p_\infty(x)$ has
arbitrarily many realizations producing models for each element in
$\mathcal{T}$. Thus there are only two possibilities for open sets
$U$: either $U=\varnothing$ or $U=\mathcal{T}$.

\medskip
{\bf Remark 2.15.} Let $\mathcal{T}_{\rm fin}$ be the class of all
theories for finite structures. By compactness, for a set
$\mathcal{T}\subset\mathcal{T}_{\rm fin}$, ${\rm
Cl}_E(\mathcal{T})$ is a subset of $\mathcal{T}_{\rm fin}$ if and
only if models of $\mathcal{T}$ have bounded cardinalities,
whereas ${\rm Cl}_P(\mathcal{T})$ is a subset of $\mathcal{T}_{\rm
fin}$ if and only if $\mathcal{T}$ is finite. Proposition 2.2 and
its $P$-analogue allows to describe both ${\rm Cl}_E(\mathcal{T})$
and ${\rm Cl}_P(\mathcal{T})$, in particular, the sets ${\rm
Cl}_E(\mathcal{T})\setminus\mathcal{T}_{\rm fin}$ and ${\rm
Cl}_P(\mathcal{T})\setminus\mathcal{T}_{\rm fin}$. Clearly, there
is a broad class of theories in $\overline{\mathcal{T}}$ which do
not lay in $$\bigcup\limits_{\mathcal{T}\subset\mathcal{T}_{\rm
fin}}{\rm
Cl}_E(\mathcal{T})\cup\bigcup\limits_{\mathcal{T}\subset\mathcal{T}_{\rm
fin}}{\rm Cl}^d_P(\mathcal{T}).$$ For instance, finitely
axiomatizable theories with infinite models can not be
approximated by theories in $\mathcal{T}_{\rm fin}$ in such way.

\section{Generating subsets of $E$-closed sets}

{\bf Definition.} Let $\mathcal{T}_0$ be a closed set in a
topological space $(\mathcal{T},\mathcal{O}_E(\mathcal{T}))$. A
subset $\mathcal{T}'_0\subseteq\mathcal{T}_0$ is said to be {\em
generating}\index{Set!generating} if $\mathcal{T}_0={\rm
Cl}_E(\mathcal{T}'_0)$. The generating set $\mathcal{T}'_0$ (for
$\mathcal{T}_0$) is {\em minimal}\index{Set!generating!minimal} if
$\mathcal{T}'_0$ does not contain proper generating subsets. A
minimal generating set $\mathcal{T}'_0$ is the {\em
least}\index{Set!generating!least} is $\mathcal{T}'_0$ is
contained in each generating set for $\mathcal{T}_0$.

\medskip
{\bf Remark 3.1.} Each set $\mathcal{T}_0$ has a generating subset
$\mathcal{T}'_0$ with a cardinality $\leq{\rm
max}\{|\Sigma|,\omega\}$, where $\Sigma$ is the union of the
languages for the theories in $\mathcal{T}_0$. Indeed, the theory
$T={\rm Th}(\mathcal{A}_E)$, whose $E$-classes are models for
theories in ${\rm Cl}_E(\mathcal{T}_0)$, has a model $\mathcal{M}$
with $|M|\leq{\rm max}\{|\Sigma|,\omega\}$. The $E$-classes of
$\mathcal{M}$ are models of theories in ${\rm
Cl}_E(\mathcal{T}_0)$ and the set of these theories is the
required generating set.

\medskip
{\bf Theorem 3.2.} {\em If $\mathcal{T}'_0$ is a generating set
for a $E$-closed set $\mathcal{T}_0$ then the following conditions
are equivalent:

$(1)$ $\mathcal{T}'_0$ is the least generating set for
$\mathcal{T}_0$;

$(2)$ $\mathcal{T}'_0$ is a minimal generating set for
$\mathcal{T}_0$;

$(3)$ any theory in $\mathcal{T}'_0$ is isolated by some set
$(\mathcal{T}'_0)_\varphi$, i.e., for any $T\in\mathcal{T}'_0$
there is $\varphi\in T$ such that
$(\mathcal{T}'_0)_\varphi=\{T\}$;

$(4)$ any theory in $\mathcal{T}'_0$ is isolated by some set
$(\mathcal{T}_0)_\varphi$, i.e., for any $T\in\mathcal{T}'_0$
there is $\varphi\in T$ such that
$(\mathcal{T}_0)_\varphi=\{T\}$.}

\medskip
{\bf\em Proof.} $(1)\Rightarrow(2)$ and $(4)\Rightarrow(3)$ are
obvious.

$(2)\Rightarrow(1)$. Assume that $\mathcal{T}'_0$ is minimal but
not least. Then there is a generating set $\mathcal{T}''_0$ such
that $\mathcal{T}'_0\setminus\mathcal{T}''_0\ne\varnothing$ and
$\mathcal{T}''_0\setminus\mathcal{T}'_0\ne\varnothing$. Take
$T\in\mathcal{T}'_0\setminus\mathcal{T}''_0$.

We assert that $T\in{\rm Cl}_E(\mathcal{T}'_0\setminus\{T\})$,
i.e., $T$ is an accumulation point of
$\mathcal{T}'_0\setminus\{T\}$. Indeed, since
$\mathcal{T}''_0\setminus\mathcal{T}'_0\ne\varnothing$ and
$\mathcal{T}''_0\subset{\rm Cl}_E(\mathcal{T}'_0)$, then by
Proposition 2.1, (3), $\mathcal{T}'_0$ is infinite and by
Proposition 2.2 it suffices to prove that for any $\varphi\in T$,
$(\mathcal{T}'_0\setminus\{T\})_\varphi$ is infinite. Assume on
contrary that for some $\varphi\in T$,
$(\mathcal{T}'_0\setminus\{T\})_\varphi$ is finite. Then
$(\mathcal{T}'_0)_\varphi$ is finite and, moreover, as
$\mathcal{T}'_0$ is generating for $\mathcal{T}_0$, by Proposition
2.2, $(\mathcal{T}_0)_\varphi$ is finite, too. So
$(\mathcal{T}''_0)_\varphi$ is finite and, again by Proposition
2.2, $T$ does not belong to ${\rm Cl}_E(\mathcal{T}''_0)$
contradicting to ${\rm Cl}_E(\mathcal{T}''_0)=\mathcal{T}_0$.

Since $T\in{\rm Cl}_E(\mathcal{T}'_0\setminus\{T\})$ and
$\mathcal{T}'_0$ is generating for $\mathcal{T}_0$, then
$\mathcal{T}'_0\setminus\{T\}$ is also generating for
$\mathcal{T}_0$ contradicting the minimality of
$\mathcal{T}'_0$.

$(2)\Rightarrow(3)$. If $\mathcal{T}'_0$ is finite then by
Proposition 2.1 (3), $\mathcal{T}'_0=\mathcal{T}_0$. Since
$\mathcal{T}_0$ is finite then for any $T\in \mathcal{T}_0$ there
is a formula $\varphi\in T$ negating all theories in
$\mathcal{T}_0\setminus\{T\}$. Therefore,
$(\mathcal{T}_0)_\varphi=(\mathcal{T}'_0)_\varphi$ is a singleton
containing $T$ and thus, $(\mathcal{T}'_0)_\varphi$ isolates $T$.

Now let $\mathcal{T}'_0$ be infinite. Assume that some
$T\in\mathcal{T}'_0$ is not isolated by the sets
$(\mathcal{T}'_0)_\varphi$. It implies that for any $\varphi\in
T$, $(\mathcal{T}'_0\setminus\{T\})_\varphi$ is infinite. Using
Proposition 2.2 we obtain $T\in{\rm
Cl}_E(\mathcal{T}'_0\setminus\{T\})$ contradicting the minimality
of $\mathcal{T}'_0$.

$(3)\Rightarrow(2)$. Assume that any theory $T$ in
$\mathcal{T}'_0$ is isolated by some set
$(\mathcal{T}'_0)_\varphi$. By Proposition 2.2 it implies that
$T\notin{\rm Cl}_E(\mathcal{T}'_0\setminus\{T\})$. Thus,
$\mathcal{T}'_0$ is a minimal generating set for $\mathcal{T}_0$.

$(3)\Rightarrow(4)$ is obvious for finite $\mathcal{T}'_0$. If
$\mathcal{T}'_0$ is infinite and any theory $T$ in
$\mathcal{T}'_0$ is isolated by some set
$(\mathcal{T}'_0)_\varphi$ then $T$ is isolated by the set
$(\mathcal{T}_0)_\varphi$, since otherwise using Proposition 2.2
and the property that $\mathcal{T}'_0$ generates $\mathcal{T}_0$,
there are infinitely many theories in $\mathcal{T}'_0$ containing
$\varphi$ contradicting $|(\mathcal{T}'_0)_\varphi|=1$.~$\Box$

\medskip
The equivalences $(2)\Leftrightarrow(3)\Leftrightarrow(4)$ in
Theorem 3.2 were noticed by E.A.~Palyutin.

\medskip
Theorem 3.2 immediately implies

\medskip
{\bf Corollary 3.3.} {\em For any structure $\mathcal{A}_E$,
$\mathcal{A}_E$ is $e$-minimal if and only if $\mathcal{A}_E$ is
$e$-least.}

\medskip
{\bf Definition.} Let $T$ be the theory ${\rm Th}(\mathcal{A}_E)$,
where $\mathcal{A}_E={\rm Comb}_E(\mathcal{A}_i)_{i\in I}$,
$\{{\rm Th}(\mathcal{A}_i)\mid i\in I\}=\mathcal{T}_0$. We say
that $T$ has a {\em minimal}/{\em least} {\em generating
set}\index{Set!generating!}\index{Set!generating!minimal}\index{Set!generating!least}
if ${\rm Cl}_E(\mathcal{T}_0)$ has a minimal/least generating set.

Since by Theorem 3.2 the notions of minimality and to be least
coincide in the context, below we shall consider least generating
sets as well as $e$-least structures in cases of minimal
generating sets.

\medskip
{\bf Proposition 3.4.} {\em For any closed nonempty set
$\mathcal{T}_0$ in a topological space
$(\mathcal{T},\mathcal{O}_E(\mathcal{T}))$ and for any
$\mathcal{T}'_0\subseteq\mathcal{T}_0$, the following conditions
are equivalent:

$(1)$ $\mathcal{T}'_0$ is the least generating set for
$\mathcal{T}_0$;

$(2)$ any/some structure $\mathcal{A}_E={\rm
Comb}_E(\mathcal{A}_i)_{i\in I}$, where $\{{\rm
Th}(\mathcal{A}_i)\mid i\in I\}=\mathcal{T}'_0$, is an $e$-least
model of the theory ${\rm Th}(\mathcal{A}_E)$ and $E$-classes of
each/some $e$-largest model of ${\rm Th}(\mathcal{A}_E)$ form
models of all theories in $\mathcal{T}_0$;

$(3)$ any/some structure $\mathcal{A}_E={\rm
Comb}_E(\mathcal{A}_i)_{i\in I}$, where $\{{\rm
Th}(\mathcal{A}_i)\mid i\in I\}=\mathcal{T}'_0$,
$\mathcal{A}_i\not\equiv\mathcal{A}_j$ for $i\ne j$, is an
$e$-least model of the theory ${\rm Th}(\mathcal{A}_E)$, where
$E$-classes of $\mathcal{A}_E$ form models of the least set of
theories and $E$-classes of each/some $e$-largest model of ${\rm
Th}(\mathcal{A}_E)$ form models of all theories in
$\mathcal{T}_0$.}

\medskip
{\bf\em Proof.} $(1)\Rightarrow (2)$. Let $\mathcal{T}'_0$ be the
least generating set for $\mathcal{T}_0$. Consider the structure
$\mathcal{A}_E={\rm Comb}_E(\mathcal{A}_i)_{i\in I}$, where
$\{{\rm Th}(\mathcal{A}_i)\mid i\in I\}=\mathcal{T}'_0$. Since
$\mathcal{T}'_0$ is the least generating set for $\mathcal{T}_0$,
then $\mathcal{A}_E$ is an $e$-least model of the theory ${\rm
Th}(\mathcal{A}_E)$. Moreover, by Proposition 2.2, $E$-classes of
models of ${\rm Th}(\mathcal{A}_E)$ form models of all theories in
$\mathcal{T}_0$. Thus, $E$-classes of $\mathcal{A}_E$ form models
of the least set $\mathcal{T}'_0$ of theories such that
$E$-classes of each/some $e$-largest model of ${\rm
Th}(\mathcal{A}_E)$ form models of all theories in
$\mathcal{T}_0$.

Similarly, constructing $\mathcal{A}_E$ with
$\mathcal{A}_i\not\equiv\mathcal{A}_j$ for $i\ne j$, we obtain
$(1)\Rightarrow (3)$.

Since $(3)$ is a particular case of $(2)$, we have $(2)\Rightarrow
(3)$.

$(3)\Rightarrow (1)$. Let $\mathcal{A}_E$ be an $e$-least model of
the theory ${\rm Th}(\mathcal{A}_E)$ and $E$-classes of each/some
$e$-largest model of ${\rm Th}(\mathcal{A}_E)$ form models of all
theories in $\mathcal{T}_0$. Then by the definition of ${\rm
Cl}_E$, $\mathcal{T}'_0$ is the least generating set for
$\mathcal{T}_0$.~$\Box$

\medskip
Note that any prime structure $\mathcal{A}_E$ (or a structure with
finitely many $E$-classes, or a prime structure extended by
finitely many $E$-classes), is $e$-minimal forming, by its
$E$-classes, the least generating set $\mathcal{T}'_0$ of theories
for the set $\mathcal{T}_0$ of theories corresponding to
$E$-classes of $e$-largest $\mathcal{A}'_E\equiv\mathcal{A}_E$.
Indeed, if a set $\mathcal{T}''_0$ is generating for
$\mathcal{T}_0$ then by Proposition 2.2 there is a model
$\mathcal{M}$ of $T$ consisting of $E$-classes with the set of
models such that their theories form the set $\mathcal{T}''_0$.
Since $\mathcal{A}_E$ prime (or with finitely many $E$-classes, or
a prime structure extended by finitely many $E$-classes), then
$\mathcal{A}_E$ is elementary embeddable into $\mathcal{M}$
(respectively, has $E$-classes with theories forming
$\mathcal{T}''_0$, or elementary embeddable to a restriction
without finitely many $E$-classes), then
$\mathcal{T}'_0\subseteq\mathcal{T}''_0$, and so $\mathcal{T}'_0$
is the least generating set for $\mathcal{T}_0$. Thus, Proposition
3.4 implies

\medskip
{\bf Corollary 3.5.} {\em Any theory ${\rm Th}(\mathcal{A}_E)$
with a prime model $\mathcal{M}$, or with a finite set $\{{\rm
Th}(\mathcal{A}_i)\mid i\in I\}$, or both with $E$-classes for
$\mathcal{M}$ and $\mathcal{A}_i$, has the least generating set.}

\medskip
Clearly, the converse for prime models does not hold, since finite
sets $\mathcal{T}_0$ are least generating whereas theories in
$\mathcal{T}_0$ can be arbitrary, in particular, without prime
models. Again the converse for finite sets does not hold since
there are prime models with infinite $\mathcal{T}_0$. Finally the
general converse is not true since we can combine a theory $T$
having a prime model with infinite $\mathcal{T}_0$ and a theory
$T'$ with infinitely many $E$-classes of disjoint languages and
without prime models for these classes. Denoting by
$\mathcal{T}'_0$ the set of theories for these $E$-classes, we get
the least infinite generating set
$\mathcal{T}_0\cup\mathcal{T}'_0$ for the combination of $T$ and
$T'$, which does not have a prime model.

Replacing $E$-combinations by $P$-combinations we obtain the
notions of (minimal/least) generating
set\index{Set!generating}\index{Set!generating!minimal}\index{Set!generating!least}
for ${\rm Cl}_P(\mathcal{T}_0)$.

The following example shows that Corollary 3.5 does not hold even
for disjoint $P$-combinations.

\medskip
{\bf Example 3.6.} Take structures $\mathcal{A}_{i}$,
$i\in(\omega+1)\setminus\{0\}$, in Remark 2.13 and the theories
$T_i={\rm Th}(\mathcal{A}_i)$ forming the ${\rm Cl}^{d}_P$-closed
set $\mathcal{T}$. Since $\mathcal{T}$ is generated by any its
infinite subset, we get that having prime models of ${\rm
Th}(\mathcal{A}_P)$, the closure ${\rm Cl}^{d}_P(\mathcal{T})$
does not have minimal generating sets.

For the example above, with the empty language, ${\rm
Cl}^{d,r}_P(\mathcal{T})$ is generated by any singleton $\{T\}\in
{\rm Cl}^{d,r}_P(\mathcal{T})$ since the type $p_\infty(x)$ has
arbitrarily many realizations producing models for each $T_i$,
$i\in(\omega+1)\setminus\{0\}$. Thus, each element of ${\rm
Cl}^{d,r}_P(\mathcal{T})$ forms a minimal generating set.

Adding to the language $\Sigma$ countably many unary predicate
symbols $R^{(1)}_i$, $i\in\omega\setminus\{0\}$, for constants and
putting each singleton $R_i$ into $\mathcal{A}_{i}$,
$i\in\omega\setminus\{0\}$, we get examples of ${\rm
Cl}^{d}_P(\mathcal{T})$ and ${\rm Cl}^{d,r}_P(\mathcal{T})$ with
the least (infinite) generating sets. Thus, the property of the
non-existence of minimal/least generating sets is not preserved
under expansions of theories.

We again obtain the non-existence of minimal/least generating sets
for ${\rm Cl}^{d}_P$ and ${\rm Cl}^{d,r}_P$, respectively,
expanding theories $T_i$ in the previous example by singletons
$R_j$, $j\ne i$, which are equal to $R_i$.

\medskip
Natural questions arise concerning minimal generating sets:

\medskip
{\bf Question 1.} {\em What is a characterization for the
existence of least generating sets?}

\medskip
{\bf Question 2.} {\em Is there exists a theory ${\rm
Th}(\mathcal{A}_E)$ without the least generating set?}

\medskip
Obviously, for $E$-combinations, Question 1 has an answer in terms
of Proposition 2.2 (clarified in Theorem 3.2) taking the least,
under inclusion, set $\mathcal{T}'_0$ generating the set ${\rm
Cl}_E(\mathcal{T}'_0)$. It means that $\mathcal{T}'_0$ does not
have accumulation points inside $\mathcal{T}'_0$ (with respect to
the sets $(\mathcal{T}'_0)_\varphi$), i.e., any element in
$\mathcal{T}'_0$ is isolated by some formula, whereas each element
$T$ in ${\rm Cl}_E(\mathcal{T}'_0)\setminus\mathcal{T}'_0$ is an
accumulation point of $\mathcal{T}'_0$ (again with respect to
$(\mathcal{T}'_0)_\varphi$), i.e., $\mathcal{T}'_0$ is dense in
its $E$-closure.

Note that a positive answer to Question 2 for ${\rm Cl}_P$ is
obtained in Remark 2.13.

Below we will give a more precise formulation for this answer
related to $E$-combinations and answer Question 2 for special
cases with languages.

\section{Language uniform theories and related $E$-closures}

\medskip
{\bf Definition.} A theory $T$ in a predicate language $\Sigma$ is
called {\em language uniform},\index{Theory!language uniform} or a
{\em {\rm LU}-theory}\index{{\rm LU}-theory} if
for each arity $n$
any substitution on the set of non-empty $n$-ary predicates
preserves $T$. The {\rm LU}-theory $T$ is called {\em {\rm
IILU}-theory}\index{{\rm IILU}-theory} if it has non-empty
predicates and as soon as there is a non-empty $n$-ary predicate
then there are infinitely many non-empty $n$-ary predicates and
there are infinitely many empty $n$-ary predicates.

\medskip
Below we point out some basic examples of {\rm LU}-theories:

\medskip
${\small\bullet}$ Any theory $T_0$ of infinitely many independent
unary predicates $R_k$ is a {\rm LU}-theory; expanding $T_0$ by
infinitely many empty predicates $R_l$ we get a {\rm IILU}-theory
$T_1$.

\medskip
${\small\bullet}$ Replacing independent predicates $R_k$ for $T_0$
and $T_1$ by disjoint unary predicates $R'_k$ with a cardinality
$\lambda\in(\omega+1)\setminus\{0\}$ such that each $R'_k$ has
$\lambda$ elements; the obtained theories are denoted by
$T^\lambda_0$ and $T^\lambda_1$ respectively; here, $T^\lambda_0$
and $T^\lambda_1$ are {\rm LU}-theories, and, moreover,
$T^\lambda_1$ is a {\rm IILU}-theory; we denote $T^1_0$ and
$T^1_1$ by $T^c_0$ and $T^c_1$; in this case nonempty predicates
$R'_k$ are singletons symbolizing constants which are replaced by
the predicate languages.

\medskip
${\small\bullet}$ Any theory $T$ of equal nonempty unary
predicates $R_k$ is a {\rm LU}-theory; 

\medskip
${\small\bullet}$ Similarly, {\rm LU}-theories and {\rm
IILU}-theories can be constructed using $n$-ary predicate symbols
of arbitrary arity $n$.

\medskip
${\small\bullet}$ The notion of language uniform theory can be
extended for an arbitrary language taking graphs for language
functions; for instance, theories of free algebras can be
considered as {\rm LU}-theories.

\medskip
${\small\bullet}$ Acyclic graphs with colored edges (arcs), for
which all vertices have same degree with respect to each color,
has {\rm LU}-theories. If there are infinitely many colors and
infinitely many empty binary relations then the colored graph has
a {\rm IILU}-theory.

\medskip
${\small\bullet}$ Generic arc-colored graphs without colors for
vertices \cite{SuCCMCT, Su043}, free polygonometries of free
groups \cite{SuGP}, and cube graphs with coordinated colorings of
edges \cite{SuMCT} have {\rm LU}-theories.

\medskip
The simplest example of a theory, which is not language uniform,
can be constructed taking two nonempty unary predicates $R_1$ and
$R_2$, where $R_1\subset R_2$. More generally, if a theory $T$,
with nonempty predicates $R_i$, $i\in I$, of a fixed arity, is
language uniform then cardinalities of $R_{i_1}^{\delta_1}(\bar
x)\wedge\ldots\wedge R_{i_j}^{\delta_1}(\bar x)$ do not depend on
pairwise distinct $i_1,\ldots,i_j$.

\medskip
{\bf Remark 4.1.} Any countable theory $T$ of a predicate language
$\Sigma$ can be transformed to a {\rm LU}-theory $T'$. Indeed,
since without loss of generality $\Sigma$ is countable consisting
of predicate symbols $R^{(k_n)}_n$, $n\in\omega$, then we can
step-by-step replace predicates $R_n$ by predicates $R'_n$ in the
following way. We put $R'_0\rightleftharpoons R_0$. If predicates
$R'_0,\ldots,R'_n$ of arities $r_0<\ldots<r_n$, respectively, are
already defined, we take for $R'_{n+1}$ a predicate of an arity
$r_{n+1}>{\rm max}\{r_n,k_{n+1}\}$, which is obtained from
$R'_{n+1}$ adding $r_{n+1}-k_{n+1}$ fictitious variables
corresponding to the formula
$$R'(x_1,\ldots,x_{k_{n+1}})\wedge(x_{k_{n+2}}\approx
x_{k_{n+2}})\wedge(x_{r_{n+1}}\approx x_{r_{n+1}}).$$

If the resulted {\rm LU}-theory $T'$ has non-empty predicates, it
can be transformed to a countable {\rm IILU}-theory $T''$ copying
these non-empty predicated with same domains countably many times
and adding countably many empty predicates for each arity $r_n$.

Clearly, the process of the transformation of $T$ to $T'$ do not
hold for for uncountable languages, whereas any {\rm LU}-theory
can be transformed to an {\rm IILU}-theory as above.

\medskip
{\bf Definition.} Recall that theories $T_0$ and $T_1$ of
languages $\Sigma_0$ and $\Sigma_1$ respectively are said to be
\emph{similar}\index{Theories!similar} if for any models ${\cal
M}_i\models T_i$, $i=0,1$, there are formulas of $T_i$, defining
in ${\cal M}_i$ predicates, functions and constants of language
$\Sigma_{1-i}$ such that the corresponding structure of
$\Sigma_{1-i}$ is a model of $T_{1-i}$.

Theories $T_0$ and $T_1$ of languages $\Sigma_0$ and $\Sigma_1$
respectively are said to be \emph{language
similar}\index{Theories!language similar} if $T_0$ can be obtained
from $T_1$ by some bijective replacement of language symbols in
$\Sigma_1$ by language symbols in $\Sigma_0$ (and vice versa).

\medskip
Clearly, any language similar theories are similar, but not vice
versa. Note also that, by the definition, any ${\rm LU}$-theory
$T$ is language similar to any theory $T^\sigma$ which is obtained
from $T$ replacing predicate symbols $R$ by $\sigma(R)$, where
$\sigma$ is a substitution on the set of predicate symbols in
$\Sigma(T)$ corresponding to nonempty predicates for $T$ as well
as a substitution on the set of predicate symbols in $\Sigma(T)$
corresponding to empty predicates for $T$. Thus we have

\medskip
{\bf Proposition 4.2.} {\em Let $T_1$ and $T_2$ be {\rm
LU}-theories of same language such that $T_2$ is obtained from
$T_1$ by a bijection $f_1$ {\rm (}respectively $f_2${\rm )}
mapping {\rm (}non{\rm )}empty predicates for $T_1$ to {\rm
(}non{\rm )}empty predicates for $T_2$. Then $T_1$ and $T_2$ are
language similar.}

\medskip
{\bf Corollary 4.3.} {\em Let $T_1$ and $T_2$ be countable {\rm
IILU}-theories of same language such that the restriction $T'_1$
of $T_1$ to non-empty predicates is language similar to the
restriction $T'_2$ of $T_2$ to non-empty predicates. Then $T_1$
and $T_2$ are language similar.}

\medskip
{\bf\em Proof.} By the hypothesis, there is a bijection $f_2$ for
non-empty predicates of $T_1$ and $T_2$. Since $T_1$ and $T_2$ be
countable {\rm IILU}-theories then $T_1$ and $T_2$ have countably
many empty predicates of each arity with non-empty predicates,
there is a bijection $f_1$ for empty predicates of $T_1$ and
$T_2$. Now Corollary is implied by Proposition 4.2.~$\Box$

\medskip
{\bf Definition.} For a theory $T$ in a predicate language
$\Sigma$, we denote by ${\rm Supp}_\Sigma(T)$\index{${\rm
Supp}_\Sigma(T)$} the {\em support}\index{Support} of $\Sigma$ for
$T$, i.~e., the set of all arities $n$ such that some $n$-ary
predicate $R$ for $T$ is not empty.

\medskip
Clearly, if $T_1$ and $T_2$ are language similar theories, in
predicate languages $\Sigma_1$ and $\Sigma_2$ respectively, then
${\rm Supp}_{\Sigma_1}(T_1)={\rm Supp}_{\Sigma_2}(T_2)$.

\medskip
{\bf Definition.} Let $T_1$ and $T_2$ be language similar theories
of same language $\Sigma$. We say that $T_2$ {\em language
dominates}\index{Theory!language dominating} $T_1$ and write
$T_1\sqsubseteq^L T_2$\index{$T_1\sqsubseteq^L T_2$} if for any
symbol $R\in\Sigma$, if $T_1\vdash\exists\bar{x}R(\bar{x})$ then
$T_2\vdash\exists\bar{x}R(\bar{x})$, i.~e., all predicates, which
are non-empty for $T_1$, are nonempty for $T_2$. If
$T_1\sqsubseteq^L T_2$ and $T_2\sqsubseteq^L T_1$, we say that
$T_1$ and $T_2$ are {\em language
domination-equivalent}\index{Theories!language
domination-equivalent} and write $T_1\sim^L T_2$.\index{$T_1\sim^L
T_2$}

\medskip
{\bf Proposition 4.4.} {\em The relation $\sqsubseteq^L$ is a
partial order on any set of ${\rm LU}$-theories.}

\medskip
{\bf\em Proof.} Since $\sqsubseteq^L$ is always reflexive and
transitive, it suffices to note that if $T_1\sqsubseteq^L T_2$ and
$T_2\sqsubseteq^L T_1$ then $T_1=T_2$. It follows as language
similar ${\rm LU}$-theories coincide having the same set of
nonempty predicates.~$\Box$

\medskip
{\bf Definition.} We say that $T_2$ {\em infinitely language
dominates}\index{Theory!infinitely language dominating} $T_1$ and
write $T_1\sqsubset^L_\infty T_2$\index{$T_1\sqsubset^L_\infty
T_2$} if $T_1\sqsubseteq^L T_2$ and for some $n$, there are
infinitely many new nonempty predicates for $T_2$ with respect to
$T_1$.

\medskip
Since there are infinitely many elements between any distinct
comparable elements in a dense order, we have

\medskip
{\bf Proposition 4.5.} {\em If a class of theories $\mathcal{T}$
has a dense order $\sqsubseteq^L$ then $T_1\sqsubset^L_\infty T_2$
for any distinct $T_1,T_2\in\mathcal{T}$ with $T_1\sqsubseteq^L
T_2$.}

\medskip
Clearly, if $T_1\sqsubseteq^L T_2$ then ${\rm
Supp}_{\Sigma}(T_1)\subseteq{\rm Supp}_{\Sigma}(T_2)$ but not vice
versa. In particular, there are theories $T_1$ and $T_2$ with
$T_1\sqsubset^L_\infty T_2$ and ${\rm Supp}_{\Sigma}(T_1)={\rm
Supp}_{\Sigma}(T_2)$.

\medskip
Let $T_0$ be a {\rm LU}-theory with infinitely many nonempty
predicate of some arity $n$, and $I_0$ be the set of indexes for
the symbols of these predicates.

Now for each infinite $I\subseteq I_0$ with $|I|=|I_0|$, we denote
by $T_I$ the theory which is obtained from the complete subtheory
of $T_0$ in the language $\{R_k\mid k\in I\}$ united with symbols
of all arities $m\ne n$ and expanded by empty predicates $R_l$ for
$l\in I_0\setminus I$, where $|I_0\setminus I|$ is equal to the
cardinality of the set empty predicates for $T_0$, of the arity
$n$.

By the definition, each $T_I$ is language similar to $T_0$: it
suffices to take a bijection $f$ between languages of $T_I$ and
$T_0$ such that (non)empty predicates of $T_I$ in the arity $n$
correspond to (non)empty predicates of $T_0$ in the arity $n$, and
$f$ is identical for predicate symbols of the arities $m\ne n$. In
particular,

Let $\mathcal{T}$ be an infinite family of theories $T_I$, and
$T_J$ be a theory of the form above (with infinite $J\subseteq
I_0$ such that $|J|=|I_0|$). The following proposition modifies
Proposition 2.2 for the $E$-closure ${\rm Cl}_E(\mathcal{T})$.

\medskip
{\bf Proposition 4.6.} {\em If $T_J\notin\mathcal{T}$ then
$T_J\in{\rm Cl}_E(\mathcal{T})$ if and only if for any finite set
$J_0\subset I_0$ there are infinitely many $T_I$ with $J\cap
J_0=I\cap J_0$.}

\medskip
{\bf\em Proof.} By the definition each theory $T_J$ is defined by
formulas describing $P_k\ne\varnothing\Leftrightarrow k\in J$.
Each such a formula $\varphi$ asserts for a finite set $J_0\subset
I_0$ that if $k\in J_0$ then $R_k\ne\varnothing\Leftrightarrow
k\in J$. It means that $\{k\in J_0\mid P_k\ne\varnothing\}=J\cap
J_0$. On the other hand, by Proposition 2.2, $T_J\in{\rm
Cl}_E(\mathcal{T})$ if and only if each such formula $\varphi$
belongs to infinitely many theories $T_I$ in $\mathcal{T}$, i.e.,
for infinitely many indexes $I$ we have $I\cap J_0=J\cap
J_0$.~$\Box$

\medskip
Now we take an infinite family $F$ of infinite indexes $I$ such
that $F$ is linearly ordered by $\subseteq$ and if $I_1\subset
I_2$ then $I_2\setminus I_1$ is infinite. The set $\{T_I\mid I\in
F\}$ is denoted by $\mathcal{T}_F$.

For any infinite $F'\subseteq F$ we denote by $\overline{\rm
lim}\,F'$\index{$\overline{\rm lim}\,F'$} the union-set $\bigcup
F'$ and by $\underline{\rm lim}\,F'$\index{$\underline{\rm
lim}\,F'$} intersection-set $\bigcap F'$. If $\overline{\rm
lim}\,F'$ (respectively $\underline{\rm lim}\,F'$) does not belong
to $F'$ then it is called the {\em upper} ({\em lower}) {\em
accumulation point}\index{Accumulation
point!upper}\index{Accumulation point!lower} (for $F'$). If $J$ is
an upper or lower accumulation point we simply say that $J$ is an
{\em accumulation point}\index{Accumulation point}.

\medskip
{\bf Corollary 4.7.} {\em If $T_J\notin\mathcal{T}_F$ then
$T_J\in{\rm Cl}_E(\mathcal{T}_F)$ if and only if $J$ is an {\rm
(}upper or lower{\rm )} accumulation point for some infinite
$F'\subseteq F$.}

\medskip
{\bf\em Proof.} If $J=\overline{\rm lim}\,F'$ or $J=\underline{\rm
lim}\,F'$ then for any finite set $J_0\subset I_0$ there are
infinitely many $T_I$ with $J\cap J_0=I\cap J_0$. Indeed, if
$J=\bigcup F'$ then for any finite $J_0\subset I_0$ there are
infinitely many $I\in F'$ such that $I\cap J_0$ contains exactly
same elements as $J\cap J_0$ since otherwise we have
$J\subset\bigcup F'$. Similarly the assertion holds for $J=\bigcap
F'$. By Proposition 4.6 we have $T_J\in{\rm Cl}_E(\mathcal{T}_F)$.

Now let $J\ne\overline{\rm lim}\,F'$ and $J\ne\underline{\rm
lim}\,F'$ for any infinite $F'\subseteq F$. In this case for each
$F'\subseteq F$,  either $J$ contains new index $j$ for a nonempty
predicate with respect to $\bigcup F'$ for each $F'\subseteq F$
with $\bigcup F'\subseteq J$ or $\bigcap F'$ contains new index
$j'$ for a nonempty predicate with respect to $J$ for each
$F'\subseteq F$ with $\bigcap F'\supseteq J$. In the first case,
for $J_0=\{j\}$ there are no $I\in F'$ such that $I\cap J_0=J\cap
J_0$. In the second case, for $J_0=\{j'\}$ there are no $I\in F'$
such that $I\cap J_0=J\cap J_0$. By Proposition 4.6 we get
$T_J\notin{\rm Cl}_E(\mathcal{T}_F)$.~$\Box$

\medskip
By Corollary 4.7 the action of the operator ${\rm Cl}_E$ for the
families $\mathcal{T}_F$ is reduced to unions and intersections of
{\em index}\index{Set!index} subsets of $F$.

Now we consider possibilities for the linearly ordered sets
$\mathcal{F}=\langle F;\subseteq\rangle$ and their closures
$\overline{\mathcal{F}}=\langle \overline{F};\subseteq\rangle$
related to ${\rm Cl}_E$.

The structure $\mathcal{F}$ is called {\em
discrete}\index{Structure!discrete} if $F$ does not contain
accumulation points.

By Corollary 4.7, if $\mathcal{F}$ is discrete then for any $J\in
F$, $T_J\notin{\rm Cl}_E(\mathcal{T}_{F\setminus\{J\}})$. Thus we
get

\medskip
{\bf Proposition 4.8.} {\em For any discrete $\mathcal{F}$,
$\mathcal{T}_F$ is the least generating set for ${\rm
Cl}_E({\mathcal{T}_F})$.}

\medskip
By Proposition 4.8, for any discrete $\mathcal{F}$,
$\mathcal{T}_F$ can be reconstructed from ${\rm
Cl}_E(\mathcal{T}_F)$ removing accumulation points, which always
exist. For instance, if $\langle F;\subseteq\rangle$ is isomorphic
to $\langle\omega;\leq\rangle$ or $\langle\omega^\ast;\leq\rangle$
(respectively, isomorphic to $\langle \mathbb Z;\leq\rangle$) then
${\rm Cl}_E(\mathcal{T}_F)$ has exactly one (two) new element(s)
$\overline{\rm lim}\,F$ or $\underline{\rm lim}\,F$ (both
$\overline{\rm lim}\,F$ and $\underline{\rm lim}\,F$).

\medskip
Consider an opposite case: with dense $\mathcal{F}$. Here, if
$\mathcal{F}$ is countable then, similarly to $\langle\mathbb
Q;\leq\rangle$, taking cuts for $\mathcal{F}$, i.~e., partitions
$(F^-,F^+)$ of $F$ with $F^- < F^+$, we get the closure
$\overline{F}$ with continuum many elements. Thus, the following
proposition holds.

\medskip
{\bf Proposition 4.9.} {\em For any dense $\mathcal{F}$,
$|\overline{F}|\geq 2^\omega$.}

\medskip
Clearly, there are dense $\mathcal{F}$ with dense and non-dense
$\overline{\mathcal{F}}$. If $\overline{\mathcal{F}}$ is dense
then, since $\overline{\overline{F}}=\overline{F}$, there are
dense $\mathcal{F}_1$ with $|F_1|=|\overline{F_1}|$. In
particular, it is followed by Dedekind theorem on completeness of
$\mathbb R$.

Answering Question 4 we have

\medskip
{\bf Proposition 4.10.} {\em If $\overline{\mathcal{F}}$ is dense
then ${\rm Cl}_E({\mathcal{T}_F})$ does not contain the least
generating set.}

\medskip
{\bf\em Proof.} Assume on contrary that ${\rm
Cl}_E({\mathcal{T}_F})$ contains the least generating set with a
set $F_0\subseteq F$ of indexes. By the minimality $F_0$ does not
contain both the least element and the greatest element. Thus
taking an arbitrary $J\in F_0$ we have that for the cut
$(F^-_{0,J},F^+_{0,J})$, where $F^-_{0,J}=\{J^-\in F_0\mid
J^-\subset J\}$ and $F^+_{0,J}=\{J^+\in F_0\mid J^+\supset J\}$,
$J=\overline{\rm lim}\,F^-_{0,J}$ and $J=\underline{\rm
lim}\,F^+_{0,J}$. Thus, $F_0\setminus\{J\}$ is again a set of
indexes for a generating set for ${\rm Cl}_E({\mathcal{T}_F})$.
Having a contradiction we obtain the required assertion.~$\Box$

\medskip
Combining Proposition 3.4 and Proposition 4.10 we obtain

\medskip
{\bf Corollary 4.11.} {\em If $\overline{\mathcal{F}}$ is dense
then ${\rm Th}(\mathcal{A}_E)$ does not have $e$-least models and,
in particular, it is not small.}

\medskip
{\bf Remark 4.12.} The condition of the density of
$\overline{\mathcal{F}}$ for Proposition 4.10 is essential.
Indeed, we can construct step-by step a countable dense structure
$\mathcal{F}$ without endpoints such that for each $J\in F$ and
for its cut $(F^-_{J},F^+_{J})$, where $F^-_{J}=\{J^-\in F\mid
J^-\subset J\}$ and $F^+_{J}=\{J^+\in F\mid J^+\supset J\}$,
$J\supset\overline{\rm lim}\,F^-_{J}$ and $J\subset\underline{\rm
lim}\,F^+_{J}$. In this case ${\rm Cl}_E({\mathcal{T}_F})$
contains the least generating set $\{T_J\mid J\in F\}$.

\medskip
In general case, if an element $J$ of $F$ has a successor $J'$ or
a predecessor $J^{-1}$ then $J$ defines a connected component with
respect to the operations $\cdot'$ and $\cdot^{-1}$. Indeed,
taking closures of elements in $F$ with respect to $\cdot'$ and
$\cdot^{-1}$ we get a partition of $F$ defining an equivalence
relation such that two elements $J_1$ and $J_2$ are equivalent if
and only if $J_2$ is obtained from $J_1$ applying $\cdot'$ or
$\cdot^{-1}$ several (maybe zero) times.

Now for any connected component $C$ we have one of the following
possibilities:

(i) $C$ is a singleton consisting of an element $J$ such that
$J\ne\overline{\rm lim}\,F^-_J$ and $J\ne\underline{\rm
lim}\,F^+_J$; in this case $J$ is not an accumulation point for
$F\setminus\{J\}$ and $T_J$ belongs to any generating set for
${\rm Cl}_E(\mathcal{T}_F)$;

(ii) $C$ is a singleton consisting of an element $J$ such that
$J=\overline{\rm lim}\,F^-_J$ or $J=\underline{\rm lim}\,F^+_J$,
and $\overline{\rm lim}\,F^-_J\ne\underline{\rm lim}\,F^+_J$; in
this case $J$ is an accumulation point for exactly one of $F^-_J$
and $F^+_J$, $J$ separates $F^-_J$ and $F^+_J$, and $T_J$ can be
removed from any generating set for ${\rm Cl}_E(\mathcal{T}_F)$
preserving the generation of ${\rm Cl}_E(\mathcal{T}_F)$; thus
$T_J$ does not belong to minimal generating sets;

(iii) $C$ is a singleton consisting of an element $J$ such that
$J=\overline{\rm lim}\,F^-_J=\underline{\rm lim}\,F^+_J$; in this
case $J$ is a (unique) accumulation point for both $F^-_J$ and
$F^+_J$, moreover, again $T_J$ can be removed from any generating
set for ${\rm Cl}_E(\mathcal{T}_F)$ preserving the generation of
${\rm Cl}_E(\mathcal{T}_F)$, and $T_J$ does not belong to minimal
generating sets;

(iv) $|C|>1$ (in this case, for any intermediate element $J$ of
$C$, $T_J$ belongs to any generating set for ${\rm
Cl}_E(\mathcal{T}_F)$), $\underline{\rm
lim}\,C\supset\overline{\rm lim}\,F^-_{\underline{\rm lim}\,C}$
and $\overline{\rm lim}\,C\subset\underline{\rm
lim}\,F^+_{\overline{\rm lim}\,C}$; in this case, for the
endpoint(s) $J^\ast$ of $C$, if it (they) exists, $T_{J^\ast}$
belongs to any generating set for ${\rm Cl}_E(\mathcal{T}_F)$;

(v) $|C|>1$, and $\underline{\rm lim}\,C=\overline{\rm
lim}\,F^-_{\underline{\rm lim}\,C}$ or $\overline{\rm
lim}\,C=\underline{\rm lim}\,F^+_{\overline{\rm lim}\,C}$; in this
case, for the endpoint $J^\ast=\underline{\rm lim}\,C$ of $C$, if
it exists, $T_{J^\ast}$ does not belong to minimal generating sets
of ${\rm Cl}_E(\mathcal{T}_F)$, and for the endpoint
$J^{\ast\ast}=\overline{\rm lim}\,C$ of $C$, if it exists,
$T_{J^{\ast\ast}}$ does not belong to minimal generating sets of
${\rm Cl}_E(\mathcal{T}_F)$.

Summarizing (i)--(v) we obtain the following assertions.

\medskip
{\bf Proposition 4.13.} {\em A partition of $F$ by the connected
components forms {\em discrete intervals} or, in particular,
singletons of $\mathcal{F}$, where only endpoints $J$ of these
intervals can be among elements $J^{\ast\ast}$ such that
$T_{J^{\ast\ast}}$ does not belong to minimal generating sets of
${\rm Cl}_E(\mathcal{T}_F)$.}

\medskip
{\bf Proposition 4.14.} {\em If $(F^-,F^+)$ is a cut of $F$ with
$\overline{\rm lim}\,F^-=\underline{\rm lim}\,F^+$ {\rm
(}respectively $\overline{\rm lim}\,F^-\subset\underline{\rm
lim}\,F^+${\rm )} then any generating set $\mathcal{T}^0$ for
${\rm Cl}_E(\mathcal{T}_F)$ is represented as a {\rm
(}disjoint{\rm )} union of generating set $\mathcal{T}^0_{F^-}$
for ${\rm Cl}_E(\mathcal{T}_{F^-})$ and of generating set
$\mathcal{T}^0_{F^+}$ for ${\rm Cl}_E(\mathcal{T}_{F^+})$,
moreover, any {\rm (}disjoint{\rm )} union of a generating set for
${\rm Cl}_E(\mathcal{T}_{F^-})$ and of a generating set for ${\rm
Cl}_E(\mathcal{T}_{F^+})$ is a generating set $\mathcal{T}^0$ for
${\rm Cl}_E(\mathcal{T}_F)$.}

\medskip
Proposition 4.14 implies

\medskip
{\bf Corollary 4.15.} {\em If $(F^-,F^+)$ is a cut of $F$ then
${\rm Cl}_E(\mathcal{T}_F)$ has the least generating set if and
only if ${\rm Cl}_E(\mathcal{T}_{F^-})$ and ${\rm
Cl}_E(\mathcal{T}_{F^+})$ have the least generating sets.}

\medskip
Considering $\subset$-ordered connected components we have that
{\em discretely ordered} intervals in $\overline{\mathcal{F}}$,
consisting of discrete connected components and their limits
$\underline{\rm lim}\,$ and $\overline{\rm lim}\,$, are alternated
with densely ordered intervals including their limits. If
$\overline{\mathcal{F}}$ contains an (infinite) dense interval,
then by Proposition 4.10, ${\rm Cl}_E(\mathcal{T}_F)$ does not
have the least generating set. Conversely, if
$\overline{\mathcal{F}}$ does not contain dense intervals then
${\rm Cl}_E(\mathcal{T}_F)$ contains the least generating set.
Thus, answering Questions 1 and 2 for ${\rm Cl}_E(\mathcal{T}_F)$,
we have

\medskip
{\bf Theorem 4.16.} {\em For any linearly ordered set
$\mathcal{F}$, the following conditions are equivalent:

$(1)$ ${\rm Cl}_E(\mathcal{T}_F)$ has the least generating set;

$(2)$ $\overline{\mathcal{F}}$ does not have dense intervals.}

\medskip
{\bf Remark 4.17.} Theorem 4.16 does not hold for some
non-linearly ordered $\mathcal{F}$. Indeed, taking countably many
disjoint copies $\mathcal{F}_q$, $q\in\mathbb Q$, of linearly
ordered sets isomorphic to $\langle\omega,\leq\rangle$ and
ordering limits $J_q=\overline{\rm lim}\,F_q$ by the ordinary
dense order on $\mathbb Q$ such that $\{J_q\mid q\in\mathbb Q\}$
is densely ordered, we obtain a dense interval $\{J_q\mid
q\in\mathbb Q\}$ whereas the set $\cup\{F_q\mid q\in\mathbb Q\}$
forms the least generating set $\mathcal{T}_0$ of theories for
${\rm Cl}_E(\mathcal{T}_0)$.

The above operation of extensions of theories for $\{J_q\mid
q\in\mathbb Q\}$ by theories for $\mathcal{F}_q$ as well as
expansions of theories of the empty language to theories for
$\{J_q\mid q\in\mathbb Q\}$ confirm that the (non)existence of a
least/minimal generating set for ${\rm Cl}_E(\mathcal{T}_0)$ is
not preserved under restrictions and expansions of theories.

\medskip
{\bf Remark 4.18.} Taking an arbitrary theory $T$ with a non-empty
predicate $R$ of an arity $n$, we can modify Theorem 4.16 in the
following way. Extending the language $\Sigma(T)$ by infinitely
many $n$-ary predicates interpreted exactly as $R$ and by
infinitely many empty $n$-ary predicates we get a class
$\mathcal{T}_{T,R}$ of theories {\em $R$-generated}\index{Class of
theories!$R$-generated} by $T$. The class $\mathcal{T}_{T,R}$
satisfies the following: any linearly ordered $\mathcal{F}$ as
above is isomorphic to some family $\mathcal{F}'$, under
inclusion, sets of indexes of non-empty predicates for theories in
$\mathcal{T}_{T,R}$ such that strict inclusions $J_1\subset J_2$
for elements in $\mathcal{F}'$ imply that cardinalities
$J_2\setminus J_1$ are infinite and do not depend on choice of
$J_1$ and $J_2$. Theorem 4.16 holds for linearly ordered
$\mathcal{F}'$ involving the given theory $T$.

\section{On $e$-spectra for families of language uniform theories}

\medskip
{\bf Remark 5.1.} Remind \cite[Proposition 4.1, (7)]{cs} that if
$T={\rm Th}(\mathcal{A}_E)$ has an $e$-least model $\mathcal{M}$
then $e$-${\rm Sp}(T)=e$-${\rm Sp}(\mathcal{M})$. Then, following
\cite[Proposition 4.1, (5)]{cs}, $e$-${\rm
Sp}(T)=|\mathcal{T}_0\setminus\mathcal{T}'_0|$, where
$\mathcal{T}'_0$ is the (least) generating set of theories for
$E$-classes of $\mathcal{M}$, and $\mathcal{T}_0$ is the closed
set of theories for $E$-classes of an $e$-largest model of $T$.
Note also that $e$-${\rm Sp}(T)$ is infinite if $\mathcal{T}_0$
does not have the least generating set.

\medskip
Remind that, as shown in \cite[Propositions 4.3]{cs}, for any
cardinality $\lambda$ there is a theory $T={\rm
Th}(\mathcal{A}_E)$ of a language $\Sigma$ such that
$|\Sigma|=|\lambda+1|$ and $e$-${\rm Sp}(T)=\lambda$. Modifying
this proposition for the class of ${\rm LU}$-theories we obtain

\medskip
{\bf Proposition 5.2.} {\em $(1)$ For any $\mu\leq\omega$ there is
an $E$-combination $T={\rm Th}(\mathcal{A}_E)$ of ${\rm
IILU}$-theories in a language $\Sigma$ of the cardinality $\omega$
such that $T$ has an $e$-least model and $e$-${\rm Sp}(T)=\mu$.

$(2)$ For any uncountable cardinality $\lambda$ there is an
$E$-combination $T={\rm Th}(\mathcal{A}_E)$ of ${\rm
IILU}$-theories in a language $\Sigma$ of the cardinality
$\lambda$ such that $T$ has an $e$-least model and $e$-${\rm
Sp}(T)=\lambda$.}

\medskip
{\bf\em Proof.} In view of Propositions 3.4, 4.8, and Remark 5.1,
it suffices to take an $E$-combination of ${\rm IILU}$-theories of
a language $\Sigma$ of the cardinality $\lambda$ and with a
discrete linearly ordered set $\mathcal{F}$ having:

1) $\mu\leq\omega$ accumulation points if $\lambda=\omega$;

2) $\lambda$ accumulation points if $\lambda>\omega$.

We get the required $\mathcal{F}$ for (1) taking:

(i) finite $F$ for $\mu=0$;

(ii) $\mu/2$ discrete connected components, forming $\mathcal{F}$,
with the ordering type $\langle\mathbb Z;\leq\rangle$ and having
pairwise distinct accumulation points, if $\mu>0$ is even natural;

(iii) $(\mu-1)/2$ discrete connected components, forming
$\mathcal{F}$, with the ordering type $\langle\mathbb
Z;\leq\rangle$ and one connected components with the ordering type
$\langle\omega;\leq\rangle$ such that all accumulation points are
distinct, if $\mu>0$ is odd natural;

(iv) $\omega$ discrete connected components, forming
$\mathcal{F}$, with the ordering type $\langle\mathbb
Z;\leq\rangle$, if $\mu=\omega$.

The required $\mathcal{F}$ for (2) is formed by (uncountably many)
$\lambda$ discrete connected components, forming $\mathcal{F}$,
with the ordering type $\langle\mathbb Z;\leq\rangle$.~$\Box$

\medskip
Combining Propositions 3.4, 4.10, Theorem 4.16, and Remark 5.1
with $\overline{\mathcal{F}}$ having dense intervals, we get

\medskip
{\bf Proposition 5.3.} {\em  For any infinite cardinality
$\lambda$ there is an $E$-combination $T={\rm Th}(\mathcal{A}_E)$
of ${\rm IILU}$-theories in a language $\Sigma$ of cardinality
$\lambda$ such that $T$ does not have $e$-least models and
$e$-${\rm Sp}(T)\geq{\rm max}\{2^\omega,\lambda\}$.}

\bigskip
\noindent Sobolev Institute of Mathematics, \\ 4, Acad. Koptyug
avenue, Novosibirsk, 630090, Russia; \\ Novosibirsk State
Technical
University, \\ 20, K.Marx avenue, Novosibirsk, 630073, Russia; \\
Novosibirsk State University, \\ 2, Pirogova street, Novosibirsk,
630090, Russia;
\\ Institute of Mathematics and Mathematical Modeling, \\
125, Pushkina Street, Almaty, 050010, Kazakhstan
\end{document}